\documentclass[preprint, a4paper, 10pt, oneside, onecolumn]{elsarticle}

\usepackage{amssymb, amsthm,amsmath}
\usepackage{graphicx}
\usepackage{subfigure}
\usepackage{float}    
\usepackage{amsthm}
\usepackage{amsmath}
\usepackage{mathrsfs}

\usepackage[misc]{ifsym}
\usepackage{diagbox}
\usepackage{graphicx} 
\usepackage{amsmath,latexsym,amssymb,amsfonts,amsbsy}
\usepackage{subfigure} 
\usepackage{algorithmic}
\usepackage{algorithm}
\usepackage{color,xcolor}
\usepackage{booktabs}
\usepackage{threeparttable}
\usepackage{multicol}
\usepackage{multirow}
\usepackage{epstopdf}
\numberwithin{equation}{section}
\newtheorem{theorem}{Theorem}[section]
\newtheorem{lemma}{Lemma}[section]

\begin{document} 

\begin{sloppypar}
    \begin{frontmatter}

        \title{Numerical analysis of a 1/2-equation model of turbulence}

        \author[1]{Wei-Wei Han}
        \ead{hanweiwei@stu.xjtu.edu.cn}
        
        \author[2]{Rui Fang} 
        \ead{rui10@pitt.edu}
       
        \author[2]{William Layton}
        \ead{wjl@pitt.edu}

    \address[1]{School of Mathematics and Statistics, Xi'an Jiaotong University, Xi'an, Shaanxi 710049, China}
    \address[2]{Department of Mathematics, Universerty of Pittsburgh, Pittsburgh, PA, 15260, USA}

        \begin{abstract}
        	The recent 1/2-equation model of turbulence is a simplification of the standard Kolmogorov-Prandtl 1-equation URANS model. Surprisingly, initial numerical tests indicated that the 1/2-equation model produces comparable velocity statistics at reduced cost. It is also a test problem and first step for developing numerical analysis to address a full 1-equation model. This report begins the numerical analysis of the 1/2 equation model. Stability, convergence and error estimates are proven for a semi-discrete and fully discrete approximation. Finally, numerical tests are conducted to validate our convergence theory.
        \end{abstract}

    \end{frontmatter}

\section{Introduction}
The numerical analysis of many aspects of laminar flows of incompressible,
viscous fluids is increasingly understood. On the other hand, almost all
fundamental issues are unresolved in the numerical analysis of turbulence
models considered useful in practice, such as URANS (Unsteady Reynolds
Averaged Navier-Stokes) models. This report gives a complete stability,
convergence and error analysis for the 1/2-equation URANS model%
\begin{gather}
	v_{t}+v\cdot\nabla v-\nabla\cdot\left(  \left[  2\nu+\mu(y)k(t)\tau\right]
	\nabla^{s}v\right)  +\nabla p=f(x)\text{ \& }\nabla\cdot v=0\text{,
	}\label{halfModel}\\
	\frac{d}{dt}k(t)+\frac{\sqrt{2}}{2}\tau^{-1}k(t)=\frac{1}{|\Omega|}%
	\int_{\Omega}\mu(y)k(t)\tau|\nabla^{s}v(x,t)|^{2}dx\text{,}\label{halfModel1}
\end{gather}
where $\nabla^{s}v=$ the symmetric part of the gradient, $\tau=$ model
time-scale, $y=$ wall normal distance to no-slip boundaries, $L=$ diameter
($\Omega$) and $\mu(y)=0.55\left(  y/L\right)  ^{2}$. This is a simplification
of the standard 1-equation model where $k(t)$ is the space average of the
1-equation approximate turbulent kinetic energy, Section 1.1 below. The
challenge in the numerical analysis is dealing with the non-monotone
nonlinearity in the highest derivative, eddy viscosity terms and with the cubic
nonlinearity in the right hand side of the $k(t)$-equation. These are all features
shared by the full 1-equation model. The numerical analysis of the full model
remains an intractable open question. Traditional methods are limited to small time or small data. For large data, turbulence occurs and develops over long times. Thus, these methods are inappropriate for handling such issues. The only previous work occurs in the important papers \cite{bernardi2002model,bernardi2004model,chacon2014numerical,rebollo2018three}, based on other model simplifications than herein.

\subsection{Model derivation}

This paper focuses on numerical analysis of algorithms for model
approximation. The detailed model derivation and a suite of tests of model accuracy are presented in \cite{fang20231}. URANS models
approximate time averages
\begin{equation}
	v(x,t)\simeq\overline{u}(x,t):=\frac{1}{\tau}\int_{t-\tau}^{t}u(x,t^{\prime
	})dt^{\prime}\label{AVelocity}%
\end{equation}
of solutions of the Navier-Stokes equations (NSE).

The most common URANS model begins with an eddy viscosity closure \cite{boussinesq1877essai,jiang2020foundations} of the time
averaged NSE with eddy viscosity given by the
Kolmogorov-Prandtl formula $\nu_{T}=\mu\sqrt{k}l$ with Kolmogorov's choice, 
$l=\sqrt{2k}\tau$, $\tau :=$ a time scale (see also \cite{kean2022clipping,kean2022prandtl,layton2018urans} for recent developments). The variable $k(x,t)$ satisfies an accepted equation
modeling the turbulent kinetic energy evolution and $l(x,t)$, the turbulent
length scale, has many different specifications of increasing complexity. We make Kolmogorov's choice, $l=\sqrt{2k}\tau$.
For the full 1-equation model \cite{davidson2015turbulence, rebollo2014mathematical, mohammadi1993analysis, pope2004ten} then results%
\begin{gather}
	v_{t}+v\cdot\nabla v-\nabla\cdot\left(  \left[  2\nu+\mu k\tau\right]
	\nabla^{s}v\right)  +\nabla q=f(x)\text{, }\nabla\cdot v=0\text{, }\nonumber\\
	k_{t}+v\cdot\nabla k-\nabla\cdot\left(  \mu k\tau\nabla k\right)  +\frac
	{\sqrt{2}}{2}\tau^{-1}k=\mu k\tau|\nabla^{s}v|^{2}\text{.}\nonumber
\end{gather}
The model studied herein (\ref{halfModel})-(\ref{halfModel1}) is obtained by space averaging
the above $k-$equation.

Error analysis of a model requires estimation of the deviation of model
solution from its discrete approximation. It thus builds on how model
uniqueness proofs estimate the deviation of model solutions. In Section 3, Theorem 3.1, we prove uniqueness of strong solutions to the model (\ref{halfModel})-(\ref{halfModel1}). The main new issue in the proof is dealing with the various nonlinearities and coupling. Section 4 proves stability, convergence and error estimates for a spatially discrete, continuous time approximation. This proof builds on the analysis of uniqueness in Section 3. Section 5 presents a fully discrete numerical analysis. Numerical tests are presented in Section 6. Since the model's accuracy was studied in \cite{fang20231}, the tests in Section 6 focus on verifying the error analysis in Section 4 and 5.

\subsection{Related work}
The main challenge herein arises from the nonlinearity in the eddy viscosity term and the right hand side of (\ref{halfModel1}). To our knowledge, the only previous large data numerical analysis of fluids model with similar non-monotone, nonlinear eddy viscosity was in \cite{bernardi2002model,bernardi2004model,chacon2014numerical}. Their work studied the models under various simplifying assumptions (different from the space averaging used to simplify herein). Our uniqueness proof in Section 3 uses a standard regularity assumption to treat the NSE convection term. This necessity reflects the fact that uniqueness of time averages $\left(\ref{AVelocity}\right)$ of solutions of the NSE is as little understood as solution uniqueness. It allows our analysis to focus on the nonlinearity introduced by the turbulence model.

Finite time averaging and ensemble averaging are the two most common approaches to URANS modelling. The k-equation (1.5) studied herein was derived in \cite{fang20231} by space averaging the standard k-equation (derived independently by Prandtl \cite{prandtl1945ueber} and Kolmogorov \cite{kolmogorov1942equations}, for more details see \cite{davidson2015turbulence, rebollo2014mathematical, mohammadi1993analysis, pope2004ten,prandtl1935mechanics,cheinet2001new}). It was based on the choice of turbulence length scale $ l = \sqrt{2k}\tau$, made by Kolmogorov \cite{kolmogorov1942equations} and mentioned as an option by Prandtl \cite{prandtl1945ueber}. The idea of 1/2 equation modeling is from Johnson and King \cite{johnson1985mathematically}, see also Wilcox \cite{wilcox1998turbulence}, Section 3.7. The idea is to take a calibration parameter that must be pre-specified and allow it to be determined by local flow conditions through solving an ordinary differential equation (ODE) (considered as `1/2 equation' in turbulence modeling). In the pioneering paper Johnson and King \cite{johnson1985mathematically} posed the ODE in the streamwise variable x. In the derivation in \cite{fang20231} and herein, the ODE is formulated in time.

\section{Preliminaries}
The common Sobolev spaces and Lebesgue spaces on $\Omega$ will be denoted by $W^{k,p}(\Omega)$ and $L^p(\Omega)$ respectively \cite{2003sobolev}, equipped with the norms $\| \cdot \|_{k,p}$ and $ \| \cdot \|_{L^p} $. As for $p=2$, we adopt $W^{k}(\Omega)$ equipped with the norm $\| \cdot \|_k$ to replace $W^{k,p}(\Omega)$ for short. With respect to the $L^2 \left(\Omega\right)$ space, $\left(\cdot,\cdot\right)$ and $\| \cdot \|$  will indicate the inner product and norm. Given a Banach space $X$, we define the norm on $L^p(0,T;X)$ 
\begin{equation*}
	\|\cdot \|_{L^p (X)} := \left( \int_{0}^{T} \|\cdot \|_{X}^{p} dt \right)  ^{1/p}  \quad \text{if} \quad 1 \leq p < \infty, \qquad  \text{ess} \sup_{t \in [0,T]} \| \cdot \| _{X} \quad \text{if} \quad p = \infty. 
\end{equation*}

We will also need discrete time notations. Let $\Delta t >0$ be the time step and $t_n = n\Delta t, n = 0,1,\cdots,N = \frac{T}{\Delta t}$. Given a Banach space $X$, we define the following norms:
\begin{equation*}
	\| v \|_{l^2 \left( X \right) } := \left(\Delta t\ \sum_{n=0}^{N}\|v^n  \|_X^2 \right)^{1/2} \quad \text{and} \quad \| v \|_{l^{\infty}\left(X\right)} := \max_{0\leq n \leq N } \| v^n \|_X .
\end{equation*}

Define the following spaces $W$ and $Y$ for velocity and pressure, respectively.
\begin{equation*}
	\begin{aligned}
		W &:= H_0^1(\Omega)^d = \{u\in H^1(\Omega)^d:u|_{\partial \Omega}=0\}, \\
		Q &:= L_0^2(\Omega)=\big\{\varphi\in L^2(\Omega):\int_\Omega \varphi dx=0\big\},\\
		V &:= \{u \in W : (\nabla \cdot u, p) = 0, \forall p \in Q \}.
	\end{aligned}
\end{equation*}
We will employ the standard skew-symmetric trilinear form:
\begin{equation}
	\begin{aligned}
		b\left(u,v,w\right) &= \left( \left(u \cdot \nabla\right)v,w\right) + \frac{1}{2}\left(\left(\nabla \cdot u\right)v,w\right)\\
		& = \frac{1}{2}\left( \left(u \cdot \nabla\right)v,w\right) - \frac{1}{2}\left( \left(u \cdot \nabla\right)w,v\right) \quad \forall u, v, w \in W
	\end{aligned}
\end{equation}
For the trilinear form, we have the following bounds \cite{girault2012finite,layton2008introduction}:
\begin{equation}
	b\left(u,v,w\right) \leq 
	\begin{cases}
		C\|u \|^{1/2}  \| \nabla u \|^{1/2} \|\nabla  w \| \| \nabla v\|,\\
		C\|\nabla u \| \|\nabla w \| \|\nabla v\|, \\
		C\| \nabla u \| \|w\|_2\|v\|.\\

	\end{cases}
\end{equation}

The spatial discretization will use the classical finite element method. Suppose $\Omega_h$ be a regular mesh of $\Omega$ and $\bar{\Omega} = \cup_{M\in \Omega_h}M$. The finite element velocity and pressure spaces $W_h$ and $Q_h$ are:
\begin{align*}
	&W_h=\{u_h\in W\cap C^0(\Omega)^d:u_h|_M\in P_{k}(M)^d, \forall M\in \Omega_h\}, \\
	&Q_h=\{p_h\in Q\cap C^0(\Omega):p_h|_M\in P_{k-1}(M), \forall M\in \Omega_h\}.
\end{align*}
where $P_k(M)$ denotes the $k$th order polynomial space on $M$ and $k\geq 2$. We assume a quasi-uniform mesh. There hold \cite{2008The,1978The} the following properties for $(W_h, Q_h)$:
\begin{equation}\label{projection:error}
	\begin{aligned}
		&\inf_{u_h\in W_h}\{\|u-u_h\|+h\|u-u_h\|_1\}\leq Ch^{k+1}|u|_{k+1} \quad \forall u\in W \cap H^{k+1}(\Omega)^2,\\
		&\inf_{p_h\in Q_h}\|q-p_h\|\leq Ch^{k}|q|_{k}\quad \forall q\in H^{k}(\Omega)\cap Q,
	\end{aligned}
\end{equation}
in which $h=:$ maximum triangle diameter in $\Omega_h$.
Furthermore, we suppose that $W_h$ and $Q_h$ satisfy the discrete inf-sup condition:
\begin{equation}\label{lbb}
	\inf_{p_h\in Q_h}\sup_{u_h\in W_h}\dfrac{(p_h,\nabla\cdot u_h)}{\|p_h\|\|u_h\|_1}\geq \beta_0 >0.
\end{equation}
where $\beta_0$ is a constant. The discretely divergence free space is:
\begin{equation}
	V_h=\{u_h\in W_h:(\nabla\cdot u_h,p_h)=0,\forall p_h\in Q_h\}.
\end{equation}
Note that the Taylor-Hood element satisfies all the above conditions with $k = 2$.

The following discrete Gronwall's inequality from \cite{1990Finite} will be used.
\begin{lemma}\label{gronwall}
	Suppose that $G,\Delta t$, and $d_n, e_n, a_n, b_n $ (for integer $n\geq0$) be nonnegative numbers such that 
	$$e_N +\Delta t \sum_{n=0}^N d_n\leq \Delta t \sum_{n=0}^{N-1}b_n e_n+\Delta t \sum_{n=0}^N a_n+G,$$
	for  $\forall N \geq 1$ and $ \forall \Delta t > 0 $, then
	$$e_N + \Delta t \sum_{n=0}^N d_n\leq \exp\bigg( \Delta t\sum_{n=0}^{N-1}b_n \bigg)\bigg(\Delta t\sum_{n=0}^N a_n+G\bigg).$$
\end{lemma}

\section{Model Uniqueness}
In this section, we prove uniqueness of strong solutions of the 1/2-equation model $\left(\ref{halfModel} \right)-\left(\ref{halfModel1} \right)$. A proof of model uniqueness gives insight into critical terms in the model numerical analysis of Section 4 and 5. In particular, we assume that the 1/2-equation model has a solution satisfying:
\begin{equation}
	\int_{0}^{T}\| \nabla v \|^4 dt < \infty.
\end{equation} 
Herein, we use the following lemma from \cite{fang20231}.
\begin{lemma}\label{lemmabound}
	Consider the 1/2-equation model $\left(\ref{halfModel} \right)-\left(\ref{halfModel1} \right)$. Recall that $k\left(t^{\star}\right) >0$. Then $k\left(t\right) > 0$ for all $t > t^{\star}$. For strong solutions, there holds the following energy equality
	\begin{equation*}
		\frac{d}{dt}\left[ \frac{1}{|\Omega|} \int_{\Omega} \frac{1}{2}| v\left(x,t\right) |^2 dx + k\left(t\right) \right] + \frac{1}{|\Omega|} \int_{\Omega} \nu |\nabla ^s v\left(x,t\right) |^2 dx + \frac{\sqrt{2}}{2}\tau ^{-1}k\left(t\right)  = \frac{1}{|\Omega|} \int_{\Omega} f \cdot v\left(x,t\right)dx.
	\end{equation*}
	The following uniform in $T$ bounds on energy and dissipation rates hold:
	\begin{equation}
		\begin{aligned}
			\frac{1}{|\Omega|} \int_{\Omega} \frac{1}{2} | v\left(x,t\right) |^2 dx & \leq C,\\
			\frac{1}{T}\int_{0}^{T}\left\{ \frac{1}{|\Omega|} \int_{\Omega} \left(\nu + \nu_T\right) |\nabla ^s v\left(x,t\right) |^2 dx \right\}dt &\leq C,\\
			\frac{1}{|\Omega|} \int_{\Omega} \frac{1}{2} | v\left(x,t\right) |^2 dx + k\left(T\right) & \leq C,\\
			\frac{1}{T}\int_{0}^{T}\left\{ \frac{1}{|\Omega|} \int_{\Omega} \nu  |\nabla ^s v\left(x,t\right) |^2 dx + \frac{\sqrt{2}}{2}\tau ^{-1}k\left(t\right) \right\}dt &\leq C,
		\end{aligned}
	\end{equation}
	where the constant $C< \infty$ depends on $f, v_0\left(x\right), k\left(t^*\right),\nu,T$.
\end{lemma}
The main result of this section is as follows. 
\begin{theorem}
	Assume that $(3.1)$ holds for a solution of $\left(\ref{halfModel} \right)-\left(\ref{halfModel1} \right)$. Then the solution is unique.
\end{theorem}
\begin{proof}
	To begin, let $\left(v_1 ,k_1\right)$ and $\left(v_2 ,k_2\right)$ be two different solutions of $\left(\ref{halfModel}\right)-\left(\ref{halfModel1}\right)$ with the same data. Set $\phi = v_1 - v_2$ and $e\left(t\right) = k_1\left(t\right) - k_2\left(t\right)$. By subtraction, there holds
	\begin{equation}\label{velocity1}
		\left\{
		\begin{aligned}
			&\frac{1}{2}\frac{d}{dt}\| \phi \|^2 + \nu \| \nabla \phi \|^2 + \int_{\Omega}\phi \cdot \nabla v_1 \cdot \phi dx + \mu \tau \int_{\Omega}\left( k_1 \nabla v_1 - k_2 \nabla v_2 \right) : \nabla \phi dx = 0,\\
			&\frac{d}{dt}e\left(t\right) + \frac{\sqrt{2}}{2}\tau ^{-1}e\left(t\right) = \varepsilon_1\left(t\right) - \varepsilon_2\left(t\right),
		\end{aligned}
		\right.
	\end{equation}
	where $\varepsilon _i \left(t\right) = \frac{1}{|\Omega|}\int_{\Omega}\mu \tau k_i  |\nabla v_i|^2 dx $ since $2\| \nabla ^s v \| = \| \nabla v \|$ from $\nabla \cdot v =0$.
	
	\noindent Let the two key model terms be denoted by
	\begin{equation*}
		A := \mu \tau \int_{\Omega}\left( k_1 \nabla v_1 - k_2 \nabla v_2 \right) : \nabla \phi dx , \quad B := \varepsilon_1\left(t\right) - \varepsilon_2\left(t\right).
	\end{equation*}
	By adding and subtracting $\int_{\Omega} k_1 \nabla v_2  : \nabla \phi dx$, we have
	\begin{equation*}
		\begin{aligned}
			A &= \mu \tau \int_{\Omega}\left( k_1 \nabla v_1 - k_1 \nabla v_2 + k_1 \nabla v_2 - k_2 \nabla v_2 \right) : \nabla \phi dx\\
			& = \mu \tau \int_{\Omega}k_1 \| \nabla \phi \|^2 dx + \mu \tau \left(k_1 - k_2\right) \int_{\Omega} \nabla v_2:\nabla \phi dx \\
			&\geq \mu \tau \int_{\Omega}k_1 \| \nabla \phi \|^2 dx - \frac{\nu}{2}\| \nabla \phi \|^2 - \frac{\mu ^2 \tau ^2}{2\nu}\| \nabla v_2 \|^2 \left(k_1 -k_2\right)^2\left(t\right).
		\end{aligned}
	\end{equation*}
	Substituting the above inequality into equation $\left(\ref{velocity1}\right)$, it yields 
	\begin{equation}\label{velocity2}
		\begin{aligned}
			\frac{1}{2}\frac{d}{dt}\| \phi \|^2 &+ \frac{\nu}{2} \| \nabla \phi \|^2 + \int_{\Omega}\phi \cdot \nabla v_1 \cdot \phi dx + \mu \tau \int_{\Omega}k_1 | \nabla \phi |^2 dx \\
			&\leq \frac{\mu ^2 \tau ^2}{2\nu}\| \nabla v_2 \|^2 \left(k_1 -k_2\right)^2\left(t\right),
		\end{aligned}
	\end{equation}
	Note that the above right term as a whole belongs to $L^1\left(0,T\right)$. Since by (3.2) $k\left(t\right) \in L^{\infty}\left(0,T\right)$ and $\| \nabla v \|^2 \in L^1\left(0,T\right)$. We now need an equation for $e\left(t\right)^2 = \left(k_1 -k_2\right)^2\left(t\right)$. We have
	\begin{equation}\label{kequation2}
		\frac{1}{2}\frac{d}{dt}e^2\left(t\right) + \frac{\sqrt{2}}{2}\tau ^{-1}e^2\left(t\right) = \left( \varepsilon_1\left(t\right) - \varepsilon_2\left(t\right) \right) \cdot e\left(t\right).
	\end{equation}
	Adding $\left(\ref{velocity2}\right)$ to $\left(\ref{kequation2}\right)$ gives
	\begin{equation}\label{total}
		\begin{aligned}
			\frac{1}{2}\frac{d}{dt}\left(\| \phi \|^2 + e^2\left(t\right) \right) &+ \frac{\nu}{2} \| \nabla \phi \|^2  + \mu \tau \int_{\Omega}k_1 | \nabla \phi |^2 dx + \frac{\sqrt{2}}{2}\tau ^{-1}e^2\left(t\right) + \int_{\Omega}\phi \cdot \nabla v_1 \cdot \phi dx \\
			&\leq  \frac{\mu ^2 \tau ^2}{2\nu}\| \nabla v_2 \|^2 \left(k_1 -k_2\right)^2\left(t\right) + \left( \varepsilon_1\left(t\right) - \varepsilon_2\left(t\right) \right) \cdot e\left(t\right).
		\end{aligned}
	\end{equation}
	We now deal with the last term on the right-hand side of $\left(\ref{total}\right)$ as follows.
	\begin{equation*}
		\begin{aligned}
			&\left( \varepsilon_1\left(t\right) - \varepsilon_2\left(t\right) \right) \cdot e\left(t\right) = \mu \tau e\left(t\right) \int_{\Omega} \left( k_1\nabla v_1 :\nabla v_1 - k_1\nabla v_1 :\nabla v_2 \right)dx \\
			&+ \mu \tau e\left(t\right) \int_{\Omega}\left( k_1\nabla v_1 :\nabla v_2 - k_2\nabla v_1 :\nabla v_2 + k_2\nabla v_1 :\nabla v_2 - k_2\nabla v_2 :\nabla v_2 \right)   dx \\
			& = \mu \tau e\left(t\right) \int_{\Omega} \left( k_1\nabla v_1 :\nabla \phi + \left(k_1 - k_2\right)\nabla v_1 :\nabla v_2 + k_2\nabla v_2 :\nabla \phi \right) dx.
		\end{aligned}
	\end{equation*}
	Then 
	\begin{equation}
		\begin{aligned}\label{key}
			& \left( \varepsilon_1\left(t\right) - \varepsilon_2\left(t\right) \right) \cdot e\left(t\right) \leq e\left(t\right) \sqrt{\mu \tau \int_{\Omega}k_1 | \nabla \phi |^2 dx} \sqrt{\mu \tau \int_{\Omega}k_1 | \nabla v_1 |^2 dx} \\
			& + \mu \tau e^2\left(t\right) \sqrt{
				\int_{\Omega} | \nabla v_1 |^2 dx} \sqrt{
				\int_{\Omega} | \nabla v_2 |^2 dx} + \mu \tau e\left(t\right) \int_{\Omega}\left[\left(k_1 - k_2\right) + k_2 \right] \nabla v_2 :\nabla \phi dx \\
			& \leq \frac{\mu \tau}{2}  \int_{\Omega}k_1 | \nabla \phi |^2 dx + \left( \frac{\mu \tau}{2}  \int_{\Omega}k_1 | \nabla v_1 |^2 dx \right) e^2\left(t\right) + \mu \tau e^2\left(t\right) \| \nabla v_1 \| \| \nabla v_2 \| \\
			& +  \left( \mu \tau \| \nabla \phi \|  \| \nabla v_2 \| \right) e^2\left(t\right) + \frac{\mu \tau}{4} \int_{\Omega}k_1 | \nabla \phi |^2 dx + \left( \mu \tau \int_{\Omega}k_1 | \nabla v_2 |^2 dx \right) e^2\left(t\right).
		\end{aligned}
	\end{equation}
	Let $a\left(t\right) = \mu \tau \left( \frac{1}{2}  \int_{\Omega}k_1 | \nabla v_1 |^2 dx + \| \nabla v_1 \| \| \nabla v_2 \|  +  \| \nabla \phi \|  \| \nabla v_2 \| + \int_{\Omega}k_1 | \nabla v_2 |^2 dx \right) + \frac{\mu ^2 \tau ^2}{2\nu}\| \nabla v_2 \|^2$. Note that $a\left(t\right) \in L^1\left(0,T\right)$ by Lemma 3.1. 
	
	There remains the standard NSE term which is bounded in a standard way as in \cite{layton2008introduction,girault2012finite}:
	\begin{equation}\label{standard}
		\int_{\Omega}\phi \cdot \nabla v_1 \cdot \phi dx  \leq C\| \phi \|^{1/2} \| \nabla \phi \|^{3/2} \| \nabla v_1 \|    \leq \frac{v}{4}\| \nabla \phi  \|^2 + C\| \nabla v_1 \|^4 \| \phi \|^2 .
	\end{equation}
	Substituting $\left(\ref{key}\right)$ and $\left(\ref{standard}\right)$ into $\left(\ref{total}\right)$, it yields
	\begin{equation}
		\begin{aligned}
			\frac{1}{2}\frac{d}{dt}\left(\| \phi \|^2 + e^2\left(t\right) \right) &+ \frac{\nu}{4} \| \nabla \phi \|^2  + \frac{\mu \tau}{4} \int_{\Omega}k_1 | \nabla \phi |^2 dx + \frac{\sqrt{2}}{2}\tau ^{-1}e^2\left(t\right) \\
			&\leq  \left( C\| \nabla v_1 \|^4  + a\left(t\right) \right) \left( \| \phi \|^2 + e^2\left(t\right) \right)  \\
		\end{aligned}
	\end{equation}
	From (3.1) we have $ \left( C\| \nabla v_1 \|^4  + a\left(t\right) \right) \in L^1\left(0,T\right)$. Uniqueness follows from Gronwall's inequality.
\end{proof}

\section{The semi-discrete approximation}
The semi-discretization scheme is as follows. Suppose $v_h\left(x,0\right)$ is the given approximation of initial condition $v_0\left(x\right)$. Find $v_h :\left[0,T\right] \rightarrow W_h$ and $q_h :\left[0,T\right] \rightarrow Q_h$ for $\forall w_h \in X_h,\forall p_h \in Q_h$ satisfying
\begin{equation}\label{semiDiscretization}
	\begin{aligned}
		\left(v_{ht},w_h\right) + \left(\nu + \mu \tau k_h\left(t\right) \right)\left( \nabla v_h, \nabla w_h\right) + b\left(v_h,v_h,w_h\right) -\left(\nabla \cdot w_h ,q_h\right)  &=  \left( f,w_h \right)\\
		\left(\nabla \cdot v_h ,p_h\right) &=0.
	\end{aligned}
\end{equation}
As for the k-equation, now we have
\begin{equation}\label{semiDiscretizationTKE}
	\frac{d k_h\left(t\right) }{dt} + \frac{\sqrt{2}}{2} \tau ^{-1} k_h\left(t\right) = \frac{1}{|\Omega|}\int_{\Omega} \mu \tau k_h\left(t\right) | \nabla v_h |^2.
\end{equation}
Next, let $e := v-v_h = v- \tilde{V} - \left(v_h - \tilde{V} \right) = \eta - \phi _h, \tilde{V}\in \left\{w_h\in W_h | \left(\nabla \cdot w_h, p_h\right) = 0,  \forall p_h \in Q_h\right\}$. We begin with presenting the stability of the semi-discrete approximation as follows.
\begin{theorem}\label{semiStability}
	Under the same assumption as in Lemma \ref{lemmabound}, then there holds the following energy inequality
	\begin{equation}
		\begin{aligned}
			\frac{d}{dt}\left\{   \frac{1}{2| \Omega |} \| v_h\left(x,t\right) \|^2 +  k_h\left(t\right) \right\} &+  \frac{\nu}{| \Omega |} \|\nabla  v_h \left(x,t\right) \|^2 + \frac{\sqrt{2}}{2}\tau ^{-1}  k_h\left(t\right) \\
			&=   \frac{1}{| \Omega |} \left( f\left(x,t\right) , v_h\left(x,t\right) \right).
		\end{aligned}	
	\end{equation}
	Furthermore, there hold the following uniform bounds on the energy and dissipation rate with respect to the $T$:
	\begin{equation}
		\begin{aligned}
			\frac{1}{2}\| v_h\left(x,T\right) \|^2 \leq C < \infty, \\
			\frac{1}{T}\int_{0}^{T}\left(\nu + \mu \tau k_h\left(t\right) \right) \|\nabla  v_h \left(x,t\right) \|^2 dt \leq C < \infty, \\
			\frac{1}{2}\| v_h\left(x,T\right) \|^2 + |\Omega| k_h\left(T\right) \leq C < \infty, \\
			\frac{1}{T}\int_{0}^{T} \left( \nu \|\nabla  v_h \left(x,t\right) \|^2 + \frac{\sqrt{2}}{2}\tau ^{-1} |\Omega| k_h\left(t\right) \right) dt \leq C < \infty .
		\end{aligned}
	\end{equation}
\end{theorem}
\begin{proof}
	Taking $w_h = v_h$ in (\ref{semiDiscretization}) and using the skew-symmetric property of the trilinear term lead to:
	\begin{equation}
		\begin{aligned}
			\frac{d}{dt} \frac{1}{2}\| v_h\left(x,t\right) \|^2  & +    \left(\nu + \mu \tau k_h\left(t\right) \right) \|\nabla  v_h \left(x,t\right) \|^2  =   \left( f\left(x,t\right) , v_h\left(x,t\right) \right) \\
			& \leq \frac{\nu}{2} \|\nabla  v_h \left(x,t\right) \|^2 + \frac{1}{2\nu} \| f\left(x,t\right) \|_{-1}^2 .
		\end{aligned}
	\end{equation}
	Note that $k_h \left(t\right) $ is always nonnegative. Then a differential inequality implies that
	\begin{equation}
		\begin{aligned}
			\frac{1}{2}\| v_h\left(x,T\right) \|^2 \leq C < \infty, \\
			\frac{1}{T}\int_{0}^{T}\left(\nu + \mu \tau k_h\left(t\right) \right) \|\nabla  v_h \left(x,t\right) \|^2 dt \leq C < \infty.
		\end{aligned}
	\end{equation}
	Furthermore, taking $w_h = v_h$ in (\ref{semiDiscretization}) and multiplying by $|\Omega|$ on both side of (\ref{semiDiscretizationTKE}), then adding the two equations, it will yield:
	\begin{equation}
		\begin{aligned}
			&\frac{d}{dt}\left\{   \frac{1}{2}\| v_h\left(x,t\right) \|^2 + |\Omega| k_h\left(t\right) \right\}  +    \nu \|\nabla  v_h \left(x,t\right) \|^2 + \frac{\sqrt{2}}{2}\tau ^{-1} |\Omega| k_h\left(t\right)  \\
			& =   \left( f\left(x,t\right) , v_h\left(x,t\right) \right) 
			\leq \frac{\nu}{2} \|\nabla  v_h \left(x,t\right) \|^2 + \frac{1}{2\nu} \| f\left(x,t\right) \|_{-1}^2 .
		\end{aligned}
	\end{equation}
	Once again, a standard differential inequality leads to
	\begin{equation}
		\begin{aligned}
			\frac{1}{2}\| v_h\left(x,T\right) \|^2 + |\Omega| k_h\left(T\right) \leq C < \infty, \\
			\frac{1}{T}\int_{0}^{T} \left( \nu \|\nabla  v_h \left(x,t\right) \|^2 + \frac{\sqrt{2}}{2}\tau ^{-1} |\Omega| k_h\left(t\right) \right) dt \leq C < \infty .
		\end{aligned}
	\end{equation}
	
	
\end{proof}

\begin{theorem}\label{semiTheorem}
	Let $v$ be a sufficient smooth solution of the 1/2-equation model and in particular $\|\nabla v\| \in L^4\left(0,T\right)$. Then there holds 
	\begin{equation}
		\begin{aligned}
			&\sup_{t \in [0,T]} \| v-v_h \|^2 + \sup_{t \in [0,T]} \left(k\left(t\right) - k_h\left(t\right) \right)^2 + \nu \int_{0}^{T}\| \nabla v -\nabla v_h  \|^2 dt \\
			& \leq C \left( \| v_0 - v_h\left(0\right) \|^2 + \left(k\left(0\right) - k_h\left(0\right) \right)^2 \right)  \\
			& + C\inf_{w_h \in X_h, p_h\in Y_h}\left(\| \nabla\left(v-w_h\right) \|_{L^2\left(0,T,L^2\right)}^2  \right) + C\sup_{t \in [0,T]}\| v-w_h \|^2\\
			& + C\inf_{w_h \in X_h, p_h\in Y_h}\int_{0}^{T}\left( \| q-p_h \|^2 + \| \left(v-w_h\right)_t \|_{-1}^2 + \| \nabla v - \nabla w_h \|^2 \right)dt ,
		\end{aligned}
	\end{equation}
	where the constant $C>0$ depends on $v_0,\nu,f,T,\int_{0}^{T}\| \nabla v \|^4dt$.
\end{theorem}
\begin{proof}
	The weak formulation of the 1/2-equation model is
	\begin{equation*}
		\left(v_{t},w_h\right) + \left( \nu + \mu \tau k\left(t\right) \right)\left( \nabla v, \nabla w_h\right) + b\left(v,v,w_h\right) -\left(\nabla \cdot w_h ,q\right)  =  \left( f,w_h \right) \quad  \forall w_h \in X_h .
	\end{equation*}
	Subtracting $\left(\ref{semiDiscretization}\right)$ from the above and choosing $w_h \in \left\{w_h\in W_h | \left(\nabla \cdot w_h, p_h\right) = 0, \forall p_h \in Q_h\right\}$ yield
	\begin{equation*}
		\begin{aligned}
			\left(\phi_{ht},w_h\right) +  \nu \left( \nabla \phi_h, \nabla w_h\right) &= \left(\eta_{t},w_h\right) + \nu \left( \nabla \eta, \nabla w_h\right) + \mu \tau k\left( t \right)\left( \nabla v , \nabla w_h \right)  \\
			& - \mu \tau k_h\left( t \right)\left( \nabla v_h , \nabla w_h \right) +  b\left(v,v,w_h\right) - b\left(v_h,v_h,w_h\right) -\left(\nabla \cdot w_h ,q\right)  .
		\end{aligned}
	\end{equation*}
	
	After arranging the above equation and taking $w_h = \phi_h$, we have
	\begin{equation}\label{semierror}
		\begin{aligned}
			\frac{1}{2}\frac{d}{dt}\| \phi_h \|^2 &+ \nu  \| \nabla \phi_h \|^2 + \mu \tau k\left(t\right) \| \nabla \phi_h \|^2  = \left(\eta_{t},\phi_h\right) + \nu \left( \nabla \eta, \nabla \phi_h\right) \\ 
			& -\left(\nabla \cdot \phi_h ,q\right)  + b\left(\eta,v,\phi_h\right) - b\left(\phi_h,v,\phi_h\right) + b\left(v_h,\eta,\phi_h\right)\\
			& + \mu \tau k\left( t \right)\left( \nabla \eta, \nabla \phi_h \right) + \mu \tau \left( k\left( t \right) - k_h\left( t \right) \right) \left( \nabla v_h , \nabla \phi_h \right) = \sum_{i=1}^{8}T_i,\\
		\end{aligned}
	\end{equation}
	where we subtract and add the term $ \mu \tau k\left(t\right)\left(\nabla v_h, \nabla \phi_h \right)$. 
	
	We will bound each term of the right-hand side of equation $\left(\ref{semierror}\right)$ as follows. Using Cauchy Schwarz and Young's inequality, we have
	\begin{equation}\label{semierror1}
		T_1 = \left(\eta_{t},\phi_h\right) \leq \| \eta _t \|_{-1} \| \nabla \phi _h \| \leq \frac{\nu}{16}  \| \nabla \phi _h \|^2 + \frac{4}{\nu }\| \eta _t \|_{-1}^2,
	\end{equation}
	\begin{equation}\label{semierror2}
		T_2 = \nu \left( \nabla \eta, \nabla \phi_h\right) \leq \frac{\nu}{16}  \| \nabla \phi _h \|^2 +  4\nu  \| \nabla \eta \|^2,
	\end{equation}
	\begin{equation}\label{semierror3}
		T_3 = \left(\nabla \cdot \phi_h ,q\right) = \left(\nabla \cdot \phi_h ,q - p_h\right) \leq \frac{\nu}{16}  \| \nabla \phi_h \|^2 + C\frac{1}{ \nu}\| q-p_h \|^2 \quad \forall p_h \in Q_h.
	\end{equation}
	As for those trilinear terms, we obtain
	\begin{equation}\label{semierror4}
		T_4 = b\left(\eta,v,\phi_h\right) \leq C\| \nabla \eta \|  \| \nabla v \|  \| \nabla \phi_h \| \leq  \frac{\nu}{16} \| \nabla \phi_h \|^2 + \frac{C}{\nu}\| \nabla \eta \|^2   \| \nabla v \|^2 ,
	\end{equation}
	\begin{equation}\label{semierror5}
		T_5 = b\left(\phi_h,v,\phi_h\right) \leq C\| \phi_h \|^{1/2} \| \nabla \phi_h \|^{3/2} \| \nabla v \| \leq  \frac{\nu}{16} \| \nabla \phi_h \|^2 + C\left(\nu\right)\| \nabla v \|^4 \| \phi_h \|^2,
	\end{equation}
	\begin{equation}\label{semierror6}
		\begin{aligned}
			T_6 = b\left(v_h,\eta,\phi_h\right) &\leq C\| v_h \|^{1/2} \| \nabla v_h \|^{1/2} \| \nabla \eta \| \| \nabla \phi_h\| \\
			&\leq  \frac{\nu}{16} \| \nabla \phi_h \|^2 + C\frac{1}{\nu}\| v_h \| \| \nabla v_h \| \| \nabla \eta \|^2.
		\end{aligned}
	\end{equation}
	
	In the next, we bound those nonlinear eddy viscosity term as follows:
	\begin{equation}\label{semierror7}
		\begin{aligned}
			&T_7 = \mu \tau k\left( t \right)\left( \nabla \eta, \nabla \phi_h \right)  \leq \frac{\mu \tau}{4}k\left(t\right) \| \nabla \phi_h \|^2 + \mu \tau k\left(t\right) \| \nabla \eta \|^2 , \\
			&T_8 = \mu \tau \left( k\left( t \right) - k_h\left( t \right) \right) \left( \nabla v_h , \nabla \phi_h \right) \leq \frac{\nu}{16} \| \nabla \phi_h \|^2 + \frac{4\mu ^2 \tau ^2}{\nu} \| \nabla v_h \|^2 \left( k\left( t \right) - k_h\left( t \right) \right) ^2 .
		\end{aligned}
	\end{equation}
	
	Substituting $\left(\ref{semierror1}\right)-\left(\ref{semierror7}\right)$ into $\left(\ref{semierror}\right)$, then we obtain
	\begin{equation}\label{semierror8}
		\begin{aligned}
			\frac{1}{2}\frac{d}{dt}\| \phi_h \|^2 &+ \frac{\nu}{2} \| \nabla \phi_h \|^2 + \frac{3\mu \tau}{4} k\left(t\right) \| \nabla \phi_h \|^2 \\
			&\leq C\left(\nu\right)\| \nabla v \|^4 \| \phi_h \|^2  + \frac{C}{\nu }\| \eta _t \|_{-1}^2 +  C\left(\nu + \mu \tau k\left(t\right)\right) \| \nabla \eta \|^2  \\
			&+ \frac{C}{ \nu}\| q-p_h \|^2  + \frac{C}{\nu}  \| \nabla v \|^2 \| \nabla \eta \|^2  + \frac{C}{\nu}\| v_h \| \| \nabla v_h \| \| \nabla \eta \|^2 \\
			& + \frac{4\mu ^2 \tau ^2}{\nu} \| \nabla v_h \|^2 \left( k\left( t \right) - k_h\left( t \right) \right) ^2 .
		\end{aligned}
	\end{equation}
	
	To deal with the last term in the above equation, we need introduce the k-equation. Subtracting $\left( \ref{semiDiscretizationTKE} \right)$ from $\left( \ref{halfModel1} \right)$ and multiplying by $e\left(t\right) = \left( k\left(t\right) - k_h\left(t\right) \right)$, it yields:
	\begin{equation}\label{semierror9}
		\begin{aligned}
			\frac{1}{2}\frac{d }{dt} e\left(t\right)^2  &+ \frac{\sqrt{2}}{2} \tau ^{-1} e\left(t\right)^2 \\
			& = \mu \tau e\left(t\right) \left[ \left( k\left(t\right) \nabla v ,\nabla v  \right)   - \left( k_h\left(t\right) \nabla v_h ,\nabla v_h \right)  \right]  \\
			& = \mu \tau e\left(t\right) \left[ \left( k\left(t\right) \nabla v , \nabla v - \nabla v_h \right)  + \left(e\left(t\right) \nabla v , \nabla v_h \right)   + \left(  \nabla v - \nabla v_h ,k_h\left(t\right) \nabla v_h \right)  \right]   \\
			&  = \sum_{i=9}^{11}T_i,
		\end{aligned}
	\end{equation}
	where we add and subtract the term $ \left( k\left(t\right) \nabla v , \nabla v_h \right) $ and $ \left(  k_h\left(t\right)\nabla v  , \nabla v_h \right) $.
	
	We will bound the three terms $T_i,i=9,10,11$ as follows. As for the first term $T_9$, we have
	\begin{equation}\label{semierror10}
		\begin{aligned}
			&T_9 = \mu \tau e\left(t\right) \left( k\left(t\right) \nabla v , \nabla v - \nabla v_h \right) = \mu \tau e\left(t\right) \left( k\left(t\right) \nabla v , \nabla \eta - \nabla \phi_h \right)  \\ 
			&=  \mu \tau e\left(t\right) \left( k\left(t\right) \nabla v , \nabla \eta \right) - \mu \tau e\left(t\right) \left( k\left(t\right) \nabla v ,  \nabla \phi_h \right)\\
			& \leq \frac{\sqrt{2}}{8}\tau ^{-1} e\left(t\right)^2 + C\tau ^{3}\mu ^2 k\left(t\right)^2 \| \nabla v  \|^2 \| \nabla \eta \|^2 + \mu \tau k\left(t\right)e\left(t\right)\|  \nabla v \|  \|  \nabla \phi_h \| \\
			& \leq \frac{\sqrt{2}}{8}\tau ^{-1} e\left(t\right)^2 + C\tau ^{3}\mu ^2 k\left(t\right)^2 \| \nabla v  \|^2 \| \nabla \eta \|^2 + \frac{\mu \tau k\left(t\right)}{4} \| \nabla \phi_h \|^2 + \mu \tau k\left(t\right)\| \nabla v \|^2 e\left(t\right)^2 .
		\end{aligned}
	\end{equation}
	In a similar way, we have
	\begin{equation}\label{semierror11}
		T_{10}  =  \mu \tau e\left(t\right) \left(e\left(t\right) \nabla v , \nabla v_h \right) \leq \mu \tau \| \nabla v \|  \| \nabla v_h \|  e\left(t\right)^2 .
	\end{equation}
	and
	\begin{equation}\label{semierror12}
		\begin{aligned}
			T_{11}  &= \mu \tau e\left(t\right) \left(  \nabla v - \nabla v_h ,k_h\left(t\right) \nabla v_h \right) = \mu \tau e\left(t\right) \left(  \nabla v - \nabla v_h , \left( k\left(t\right) - e\left(t\right) \right) \nabla v_h \right) \\
			& = \mu \tau e\left(t\right) \left[  k\left(t\right) \left(  \nabla \eta  , \nabla v_h \right) - k\left(t\right) \left(  \nabla \phi_h  , \nabla v_h \right) - e\left(t\right) \left(  \nabla \eta  , \nabla v_h \right)  + e\left(t\right) \left(  \nabla \phi_h  , \nabla v_h \right) \right]   \\
			& \leq \frac{\sqrt{2}}{8}\tau ^{-1} e\left(t\right)^2 + C\tau ^{3}\mu ^2 k\left(t\right)^2 \| \nabla v_h \|^2 \| \nabla \eta \|^2 + \frac{\mu \tau k\left(t\right) }{4}  \| \nabla \phi_h \|^2 \\
			&+ C\mu \tau k\left(t\right) \| \nabla v_h \|^2 e\left(t\right)^2  + \mu \tau  \left( \| \nabla v \| + \| \nabla v_h \| \right) \| \nabla v_h \| e\left(t\right)^2.
		\end{aligned}
	\end{equation}
	Substituting $\left(\ref{semierror10}\right)-\left(\ref{semierror12}\right)$ into $\left(\ref{semierror9}\right)$ and adding $\left(\ref{semierror8}\right)$, then we obtain
	\begin{equation}
		\begin{aligned}
			& \frac{1}{2}\frac{d}{dt} \| \phi_h \|^2   + \frac{1}{2}\frac{d }{dt} e\left(t\right)^2 + \frac{\nu}{2} \| \nabla \phi_h \|^2 + \frac{\mu \tau}{4} k\left(t\right) \| \nabla \phi_h \|^2   + \frac{\sqrt{2}}{4} \tau ^{-1} e\left(t\right)^2 \\
			&\leq C\left(\nu\right)\| \nabla v \|^4 \| \phi_h \|^2 + \left( \frac{4\mu ^2 \tau ^2}{\nu} \| \nabla v_h \|^2 + \mu \tau \left( 1 + k\left(t\right)  \right)\left( \| \nabla v \|^2 + \| \nabla v_h \|^2 \right)   \right) e\left(t\right)^2 \\
			& + \frac{C}{\nu }\| \eta _t \|_{-1}^2 +  C\left(\nu + \mu \tau k\left(t\right)\right) \| \nabla \eta \|^2 + \frac{C}{ \nu}\| q-p_h \|^2  + \frac{C}{\nu}  \| \nabla v \|^2 \| \nabla \eta \|^2  \\
			&+ \frac{C}{\nu}\| v_h \| \| \nabla v_h \| \| \nabla \eta \|^2  + C\tau ^{3}\mu ^2 k\left(t\right)^2 \| \nabla v  \|^2 \| \nabla \eta \|^2 + C\tau ^{3}\mu ^2 k\left(t\right)^2 \| \nabla v_h \|^2 \| \nabla \eta \|^2  .
		\end{aligned}
	\end{equation}
	Let $b\left(t\right) = \max \left\{ C\left(\nu\right)\| \nabla v \|^4 , C\left( \frac{4\mu ^2 \tau ^2}{\nu} \| \nabla v_h \|^2 + \mu \tau \left( 1 + k\left(t\right)  \right)\left( \| \nabla v \|^2 + \| \nabla v_h \|^2 \right)   \right) \right\} $. Since $\| \nabla v \| \in L^4\left(0,T\right)$, we know $b\left(t\right) \in L^1\left(0,T\right)$, then
	\begin{equation*}
		B\left(t\right) :=\int_{0}^{t}b\left(t'\right) dt' < \infty .
	\end{equation*}
	Multiplying by the integrating factor $e^{-B\left(t\right)}$ gives
	\begin{equation}
		\begin{aligned}
			&\frac{d}{dt}\left[e^{-B\left(t\right)} \left( \| \phi_h \|^2 + e\left(t\right)^2 \right) \right] + e^{-B\left(t\right)} \nu \| \nabla \phi_h \|^2 \\
			& \leq e^{-B\left(t\right)} C \left( \| \eta _t \|_{-1}^2 +  \left(\nu + \mu \tau k\left(t\right)\right) \| \nabla \eta \|^2 + \| q-p_h \|^2  +  \| \nabla v \|^2 \| \nabla \eta \|^2  \right) \\
			&+ e^{-B\left(t\right)}C \left\{ \| v_h \| \| \nabla v_h \|   + \tau ^{3}\mu ^2 k\left(t\right)^2 \| \nabla v  \|^2 + \tau ^{3}\mu ^2 k\left(t\right)^2 \| \nabla v_h \|^2  \right\} \| \nabla \eta \|^2 .
		\end{aligned}
	\end{equation}
	
	Integrating over the time interval $\left[0,T\right]$ and multiplying by $e^{B\left(t\right)}$, it yields
	\begin{equation}
		\begin{aligned}
			&\| \phi_h \left(T\right) \|^2 + e\left(T\right)^2 + \nu \int_{0}^{T} \| \nabla \phi_h \|^2 dt \leq C\left(T,\nu\right)\| \phi_h \left(0\right) \|^2 \\
			&+ C\left(T,\nu\right)\int_{0}^{T}\left\{ \| \eta _t \|_{-1}^2 + \| q-p_h \|^2  + \| \nabla v \|^2 \| \nabla \eta \|^2     + \| v_h \| \| \nabla v_h \| \| \nabla \eta \|^2  + \| \nabla \eta \|^2  \right\}.
		\end{aligned}
	\end{equation}
	In the final, using the triangle inequality gives the final result. 
\end{proof}

\section{The fully discrete approximation}
We analyze the following full discretization scheme based on the Backward Euler (BE) time discretization. We choose the simple BE time discretization for the analysis so we can focus the analysis on the new terms in the model. Given $v_h^n, k_h^n$, finding $v_h^{n+1}\in W_h, q_h^{n+1}\in Q_h, k_h^{n+1}$ satisfies for $\forall w_h \in W_h, \forall p_h \in Q_h$
\begin{equation}\label{halferr}
	\left\{
	\begin{aligned}
		\left(\frac{v_h^{n+1}-v_h^n}{\Delta t},w_h\right) + \left(\nu + \mu \tau k_h^n\right)\left( \nabla v_h^{n+1}, \nabla w_h\right) &+ b\left(v_h^{n},v_h^{n+1},w_h\right) -\left(\nabla \cdot w_h ,q_h^{n+1}\right)  \\
		&=  \left( f^{n+1},w_h \right), \\
		\left(\nabla \cdot v_h^{n+1} ,p_h\right) &=0, \\
		\frac{k_h^{n+1}-k_h^n}{\Delta t} + \frac{\sqrt{2}}{2}\tau ^{-1}k_h^{n+1} &= \frac{1}{|\Omega|}\int_{\Omega}\mu \tau k_h^n | \nabla  v_h^{n+1}  |^2 dx.
	\end{aligned}
	\right.
\end{equation}

The main challenge of the numerical analysis is dealing with the non-monotone nonlinearity in the eddy viscosity and the right hand side of the k-equation. To streamline the analysis, we assume that the k-equation (a linear constant coefficient ODE) is solved exactly. However, it still contains errors since its right-hand side depends on the highly nonlinear energy dissipation rate and the approximation velocity. Further, $k_h\left(t\right)$  is lagged in $\nu_T$. Based on this assumption, we can obtain 
\begin{equation}\label{halferror}
	\begin{aligned}
		\left(\frac{v_h^{n+1}-v_h^n}{\Delta t},w_h\right) + \left(\nu + \mu \tau k_h\left(t_n\right)\right)\left( \nabla v_h^{n+1}, \nabla w_h\right) + b\left(v_h^{n},v_h^{n+1},w_h\right) \\
		-\left(\nabla \cdot w_h ,q_h^{n+1}\right) =  \left( f^{n+1},w_h \right) .
	\end{aligned}
\end{equation}
Before performing the derivation, let's introduce some notations. 
\begin{equation*}
	e^{n+1} := v^{n+1} - v_h^{n+1} = \left( v^{n+1} - U_h^{n+1} \right) - \left(v_h^{n+1} -  U_h^{n+1}  \right) = \eta ^{n+1} - \phi _h^{n+1} ,\quad \text{where}  \quad U_h^{n+1} \in V_h.
\end{equation*}
Taking $t = t_{n+1}$ in $\left(\ref{halfModel}\right)$ yields
\begin{equation}\label{true}
	\begin{aligned}
		\left(\frac{v^{n+1}-v^n}{\Delta t},w_h\right) &+ \left(\nu + \mu \tau k\left(t_{n+1}\right)\right)\left( \nabla v^{n+1}, \nabla w_h\right) + b\left(v^{n+1},v^{n+1},w_h\right) \\
		&-\left(\nabla \cdot w_h ,q^{n+1}\right) =  \left( f^{n+1},w_h \right) + \left(R_v^{n+1},w_h\right),
	\end{aligned}
\end{equation}
where $R_v^{n+1} = -\frac{1}{\Delta t}\int_{t^n}^{t^{n+1}}\left(t -t_{n+1}\right)v_{tt}dt.$

Before presenting the error estimates, we first show that the numerical scheme (\ref{halferror}) is unconditionally stable.

\begin{theorem}
	The scheme (\ref{halferr}) is unconditionally stable:
	\begin{equation}
		\begin{aligned}
			\| v_h^{N} \|^2 + 2| \Omega | {k_h^{n+1}} &+ \Delta t\sum_{n=0}^{N-1} \left( \nu \| \nabla v_h^{n+1} \|^2 + \sqrt{2}\tau ^{-1} |\Omega| {k_h^{n+1}} \right)  + \sum_{n=0}^{N-1} \| v_h^{n+1} - v_h^{n} \|^2  \\
			&\leq \sum_{n=0}^{N-1} C\Delta t \|f^{n+1}\|^2 + \| v_h^{0} \|^2 + 2| \Omega | {k_h^{0}} .
		\end{aligned}
	\end{equation}
\end{theorem}

\begin{proof}
	At first, taking $w_h = 2\Delta t v_h^{n+1}$ in the first equation of (\ref{halferr}), it yields
	\begin{equation}\label{stable1}
		\| v_h^{n+1} \|^2 - \| v_h^{n} \|^2 + \| v_h^{n+1} - v_h^{n} \|^2 + 2\left( \nu + \mu \tau k_h^n \right)\Delta t\| \nabla v_h^{n+1} \|^2 = 2\Delta t \left( f^{n+1}, v_h^{n+1} \right).
	\end{equation}
	Then multiply both sides of the third equation of (\ref{halferr}) by $2|\Omega|\Delta t$, we will have 
	\begin{equation}\label{stable2}
		2 | \Omega | \left( {k_h^{n+1}} -{k_h^{n}} \right) +\sqrt{2}\tau ^{-1} |\Omega| \Delta t {k_h^{n+1}} = 2\Delta t \mu \tau k_h^n \| \nabla  v_h^{n+1} \|^2 .
	\end{equation}
	Adding (\ref{stable1}) and (\ref{stable2}) and using the Young's inequality give
	\begin{equation}
		\begin{aligned}
			&\| v_h^{n+1} \|^2 - \| v_h^{n} \|^2 + \| v_h^{n+1} - v_h^{n} \|^2 + 2 | \Omega | \left( {k_h^{n+1}} -{k_h^{n}} \right) \\
			&+\sqrt{2}\tau ^{-1} |\Omega| \Delta t {k_h^{n+1}} +  \nu \Delta t\| \nabla v_h^{n+1} \|^2 \leq C\Delta t \|f^{n+1}\|^2.
		\end{aligned}
	\end{equation}
	Summing up from $n=0$ to $N-1$, we get
	\begin{equation*}
		\begin{aligned}
			\| v_h^{N} \|^2 + 2| \Omega | {k_h^{n+1}} &+ \Delta t\sum_{n=0}^{N-1} \left( \nu \| \nabla v_h^{n+1} \|^2 + \sqrt{2}\tau ^{-1} |\Omega| {k_h^{n+1}} \right)  + \sum_{n=0}^{N-1} \| v_h^{n+1} - v_h^{n} \|^2  \\
			&\leq \sum_{n=0}^{N-1} C\Delta t \|f^{n+1}\|^2 + \| v_h^{0} \|^2 + | 2\Omega | {k_h^{0}} .
		\end{aligned}
	\end{equation*}

\end{proof}

\begin{theorem}
	Suppose that the solution of 1/2-equation model is sufficiently smooth and in particular $v \in L^{\infty}\left(0,T;H^1\left(\Omega\right)\right) \cap L^2\left(0,T;H^{k+1}\left(\Omega\right)\right) \cap H^{2}\left(0,T;H^{-1}\left(\Omega\right)\right)$. Then there exists a positive constant $C\left(v,q,k,\nu,T\right)$ such that
	\begin{equation}
		\begin{aligned}
			&\| v^N - v_h^N \|^2 + \nu\Delta t\sum_{n=0}^{N-1}\| \nabla\left(v^{n+1} - v_h^{n+1}\right) \|^2 \\
			&\leq C\Delta t^2\left(  \max_{1\leq n \leq N } \| \nabla v^{n} \|^2 \left( \| k_t \|^2_{L^2\left(0,T\right)} + \| v_t \|^2_{L^2\left(0,T;H^1\right)} \right) + \| v_{tt} \|^2_{L^2\left(0,T;H^{-1}\right)}  \right) \\
			&+ Ch^{2k}\left( |\| v \| |^2_{2,k+1} + |\|  q \| |^2_{2,k} \right).
		\end{aligned}
	\end{equation}
\end{theorem}
\begin{proof}
	The analysis requires bounding several terms common to the NSE and several new terms, the analytical contribution here. The eddy viscosity terms are treated in equations $\left(\ref{totalerror5}\right)$ to $\left(\ref{totalerror7}\right)$. The terms corresponding to the energy dissipation rate errors (the right-hand side of the k-equation) are in equations $\left(\ref{totalerror7}\right)$. Subtracting $\left(\ref{halferror}\right)$ from $\left(\ref{true}\right)$ and taking $w_h \in V_h$, it yields
	\begin{equation*}
		\begin{aligned}
			\left(\frac{e^{n+1}-e^n}{\Delta t},w_h\right) &+ \nu\left(\nabla e,\nabla w_h\right) + \mu \tau k\left(t_{n+1}\right) \left(\nabla v^{n+1},\nabla w_h\right) - \mu \tau k_h\left(t_{n}\right) \left(\nabla v_h^{n+1},\nabla w_h\right) \\
			&+b\left(v^{n+1},v^{n+1},w_h\right) - b\left(v_h^{n},v_h^{n+1},w_h\right) -\left(\nabla \cdot w_h ,q^{n+1}\right) = \left(R_v^{n+1},w_h\right).
		\end{aligned}
	\end{equation*}
	By adding and subtracting $\mu \tau k\left(t_n\right)\left(\nabla v^{n+1},\nabla w_h\right)$, $b\left(v^{n},v^{n+1},w_h\right)$, $b\left(v_h^{n},v^{n+1},w_h\right)$, and taking $w_h = \phi_h^{n+1}$ in the above equation, we have
	\begin{equation}\label{totalerror}
		\begin{aligned}
			&\left(\frac{\phi_h^{n+1}-\phi_h^n}{\Delta t},\phi_h^{n+1}\right) + \nu\left(\nabla \phi_h^{n+1},\nabla \phi_h^{n+1}\right) =\left(R_v^{n+1},\phi_h^{n+1}\right)\\ 
			& + \nu\left(\nabla \eta^{n+1},\nabla \phi_h^{n+1}\right) -\left(\nabla \cdot \phi_h^{n+1} ,q^{n+1}\right) + \mu \tau \left(k\left(t_{n+1}\right) -k\left(t_{n} \right) \right)\left(\nabla v^{n+1},\nabla \phi_h^{n+1}\right)\\
			& + \mu \tau k\left(t_n\right) \left( \nabla \left( v^{n+1} - v_h^{n+1} \right),\nabla \phi _h^{n+1} \right) + \mu \tau \left( k\left(t_n\right) - k_h\left(t_n\right) \right) \left( \nabla v_h^{n+1} ,\nabla \phi _h^{n+1} \right)   \\
			&+b\left(v^{n+1}-v^{n},v^{n+1},\phi_h^{n+1}\right)  + b\left(\eta^{n},v^{n+1},\phi_h^{n+1}\right)  \\
			&+ b\left(\phi_h^{n},v^{n+1},\phi_h^{n+1}\right) + b\left(v_h^{n},\eta^{n+1},\phi_h^{n+1}\right).
		\end{aligned}
	\end{equation}
	
	In the next, we will bound each term of the right-hand side of equation $\left(\ref{totalerror}\right)$. As for the first three terms,  by using Cauchy Schwarz and Young's inequality, we have
	\begin{equation}\label{totalerror1}
		\left(R_v^{n+1},\phi_h^{n+1}\right) \leq \| R_v^{n+1} \|_{-1} \| \nabla \phi_h^{n+1} \| \leq \epsilon \nu \| \nabla \phi_h^{n+1} \|^2 + C\| R_v^{n+1} \|_{-1}^2,
	\end{equation}
	\begin{equation}\label{totalerror2}
		\nu\left(\nabla \eta^{n+1},\nabla \phi_h^{n+1}\right) \leq \| \nabla \phi_h^{n+1} \| \| \nabla \eta^{n+1} \| \leq \epsilon \nu\| \nabla \phi_h^{n+1} \|^2  +C\| \nabla \eta^{n+1} \|^2,
	\end{equation}
	\begin{equation}\label{totalerror3}
		\begin{aligned}
			\left(\nabla \cdot \phi_h^{n+1} ,q^{n+1}\right) &= \left(\nabla \cdot \phi_h^{n+1} ,q^{n+1} -p_h\right) \leq C\| \nabla \phi_h^{n+1}\| \| q^{n+1} -p_h \| \\
			& \leq \epsilon \nu \| \nabla \phi_h^{n+1}\|^2 + C\| q^{n+1} -p_h \|^2.
		\end{aligned}
	\end{equation}
	
	Furthermore, we will deal with the three terms originated from the nonlinear eddy viscosity terms, which makes it new. We have
	\begin{equation}\label{totalerror5}
		\begin{aligned}
			&\mu \tau \left(k\left(t_{n+1}\right) -k\left(t_{n} \right) \right)\left(\nabla v^{n+1},\nabla \phi_h^{n+1}\right) \\
			&\leq 	\mu \tau \left(k\left(t_{n+1}\right) -k\left(t_{n} \right) \right) \| \nabla v^{n+1} \|  \| \nabla \phi_h^{n+1} \|\\
			& \leq \epsilon \nu \| \nabla \phi_h^{n+1} \|^2 + C\| \nabla v^{n+1} \|^2 \left(k\left(t_{n+1}\right) -k\left(t_{n} \right) \right)^2 \\
			& \leq \epsilon \nu \| \nabla \phi_h^{n+1} \|^2 + C\Delta t \| \nabla v^{n+1} \|^2\int_{t_n}^{t_{n+1}}\left(k_t\right)^2 dt .
		\end{aligned}
	\end{equation}
	\begin{equation}\label{totalerror6}
		\begin{aligned}
			&\mu \tau k\left(t_n\right) \left( \nabla \left( v^{n+1} - v_h^{n+1} \right),\nabla \phi _h^{n+1} \right) \\
			&= \mu \tau k\left(t_n\right) \left( \nabla \eta ^{n+1} ,\nabla \phi _h^{n+1} \right) - \mu \tau k\left(t_n\right) \left( \nabla \phi_h ^{n+1} ,\nabla \phi _h^{n+1} \right) \\
			& \leq \frac{\mu \tau k\left(t_{n}\right)}{4} \| \nabla \phi_h^{n+1} \|^2 + \mu \tau k\left(t_{n}\right) \| \nabla \eta ^{n+1} \|^2 - \mu \tau k\left(t_n\right) \| \nabla \phi _h^{n+1} \|^2,
		\end{aligned}
	\end{equation}
	\begin{equation}\label{totalerror7}
		\begin{aligned}
			&\mu \tau \left( k\left(t_n\right) - k_h\left(t_n\right) \right) \left( \nabla v_h^{n+1} ,\nabla \phi _h^{n+1} \right) \\
			&\leq \epsilon \nu \| \nabla \phi _h^{n+1} \|^2 + C\mu ^2 \tau ^2 \|  \nabla v_h^{n+1} \|^2 \left( k\left(t_n\right) - k_h\left(t_n\right) \right)^2.
		\end{aligned}
	\end{equation}
	As for those trilinear terms, we can obtain
	\begin{equation}\label{totalerror8}
		\begin{aligned}
			b\left(v^{n+1}-v^{n},v^{n+1},\phi_h^{n+1}\right) &\leq C \| \nabla \left(v^{n+1}-v^{n} \right)  \| \| \nabla v^{n+1} \| \| \nabla \phi_h^{n+1} \| \\
			&\leq \epsilon \nu \| \nabla \phi_h^{n+1} \|^2 + C \| \nabla \left(v^{n+1}-v^{n} \right) \|^2 \| \nabla v^{n+1} \|^2 \\
			& \leq \epsilon \nu \| \nabla \phi_h^{n+1} \|^2 + C \Delta t  \| \nabla v^{n+1} \|^2 \int_{t_n}^{t_{n+1}} \| \nabla v_t \|^2 dt,
		\end{aligned}
	\end{equation}
	\begin{equation}\label{totalerror9}
		\begin{aligned}
			b\left(\eta^{n},v^{n+1},\phi_h^{n+1}\right) &\leq C\| \nabla \eta^{n}  \| \| \nabla v^{n+1} \| \| \nabla \phi_h^{n+1} \| \\
			&\leq \epsilon \nu \| \phi_h^{n+1} \|^2 + C \| \nabla \eta^{n} \| \nabla v^{n+1} \|^2,
		\end{aligned}
	\end{equation}
	\begin{equation}\label{totalerror10}
		\begin{aligned}
			b\left(\phi_h^{n},v^{n+1},\phi_h^{n+1}\right) &\leq C\|  \phi_h^{n}  \| \| v^{n+1} \|_{2} \| \nabla \phi_h^{n+1} \| \\
			&\leq \epsilon \nu \| \phi_h^{n+1} \|^2 + C    \| v^{n+1} \|_{H^2}^2 \|  \phi_h^{n}  \|^2,
		\end{aligned} 
	\end{equation}
	\begin{equation}\label{totalerror11}
		\begin{aligned}
			b\left(v_h^{n},\eta^{n+1},\phi_h^{n+1}\right) &\leq C\| v_h^{n} \|^{1/2} \| \nabla v_h^{n} \|^{1/2} \| \nabla \eta^{n+1}  \|  \| \nabla \phi_h^{n+1} \|\\
			& \leq \epsilon \nu \| \nabla \phi_h^{n+1} \|^2 + C\| v_h^{n} \| \| \nabla v_h^{n} \| \| \nabla \eta^{n+1}  \|^2.
		\end{aligned}
	\end{equation}
	Substituting $\left(\ref{totalerror1}\right)-\left(\ref{totalerror11}\right)$ into $\left(\ref{totalerror}\right)$ and taking $\epsilon = \frac{1}{18}$, then it yields
	\begin{equation}\label{totalerror12}
		\begin{aligned}
			\frac{1}{2\Delta t}\| \phi_h^{n+1} \|^2 &- \frac{1}{2\Delta t}\| \phi_h^{n} \|^2 + \frac{1}{2\Delta t}\| \phi_h^{n+1} - \phi_h^{n}  \|^2   + \frac{\nu}{2}\| \nabla \phi_h^{n+1} \|^2 + \frac{\mu \tau k\left(t_n\right)}{2} \| \nabla \phi _h^{n+1} \|^2 \\
			& \leq C    \| v^{n+1} \|_{H^2}^2 \|  \phi_h^{n}  \|^2 + C\mu ^2 \tau ^2 \|  \nabla v_h^{n+1} \|^2 \left( k\left(t_n\right) - k_h\left(t_n\right) \right)^2 \\
			& + C\| R_v^{n+1} \|_{-1}^2 + C\| \nabla \eta^{n+1} \|^2 + C\| q^{n+1} -p_h \|^2 + C\mu \tau k\left(t_{n}\right) \| \nabla \eta^{n+1} \|^2 \\
			& + C\Delta t \| \nabla v^{n+1} \|^2\int_{t_n}^{t_{n+1}}\left(k_t\right)^2 dt  + C \Delta t  \| \nabla v^{n+1} \|^2 \int_{t_n}^{t_{n+1}} \| \nabla v_t \|^2 dt\\
			&+ C \| \nabla \eta^{n} \| \nabla v^{n+1} \|^2 + C\| v_h^{n} \| \| \nabla v_h^{n} \| \| \nabla \eta^{n+1}  \|^2 .
		\end{aligned}
	\end{equation}
	Summing from $n = 0$ to $N-1$ and multiplying $2\Delta t$, we have
	\begin{equation}\label{totalerror13}
		\begin{aligned}
			\| \phi_h^{N} \|^2  &+\sum_{n=0}^{N-1}\| \phi_h^{n+1} - \phi_h^{n}  \|^2   + \Delta t\sum_{n=0}^{N-1}\left( \nu \| \nabla \phi_h^{n+1} \|^2 + \mu \tau k\left(t_n\right) \| \nabla \phi _h^{n+1} \|^2 \right) \\
			& \leq C \Delta t \sum_{n=0}^{N-1}  \left(  \| v^{n+1} \|_{H^2}^2 \|  \phi_h^{n}  \|^2 + \| R_v^{n+1} \|_{-1}^2 + \| \nabla \eta^{n+1} \|^2 \right) 
			\\
			& + C\Delta t\sum_{n=0}^{N-1}\left(  \| q^{n+1} -p_h \|^2 + \mu \tau k\left(t_{n}\right) \| \nabla \eta^{n+1} \|^2 + \| \nabla \eta^{n} \| \nabla v^{n+1} \|^2 \right) \\
			&+ C \Delta t\sum_{n=0}^{N-1}\left(  \| v_h^{n} \| \| \nabla v_h^{n} \| \| \nabla \eta^{n+1}  \|^2 + \|  \nabla v_h^{n+1} \|^2 \left( k\left(t_n\right) - k_h\left(t_n\right) \right)^2 \right) \\
			& + C\Delta t^2 \max_{1\leq n \leq N } \| \nabla v^{n} \|^2 \left( \int_{0}^{T}\left(k_t\right)^2 dt  +  \int_{0}^{T} \| \nabla v_t \|^2 dt \right).
		\end{aligned}
	\end{equation}
	
	In the final, by using the discrete Gronwall's inequality, triangle inequality and Theorem \ref{semiTheorem}, we can derive the final result.
	
\end{proof}

\section{Numerical Tests}

We use a test problem from \cite{grossmann2016high}. Consider the flow between offset circles to test the convergence rates of the numerical schemes. Due to not knowing the analytical solution, we will use the numerical results on the finer mesh as the reference solution to compute the convergence rates. The domain is a disk with a smaller offset obstacle inside. $$\Omega = \{(x,y): x^2 + y^2 \leq r_1 ^2 \cap (x-c_1)^2 + (y-c_2)^2 \geq r_2^2 \} $$ where  $r_1 =1,\quad r_2 = 0.1,\quad c= (c_1,c_2)=(\frac{1}{2},0)$. The flow is driven by a counterclockwise force $f(x,y,t) = (4x\min \left(t,1\right)(1-x^2 -y^2) , -4y\min \left(t,1\right) (1-x^2-y^2))$. Impose the no-slip boundary conditions on both circles. Herein we choose $\tau = 0.1$, $\mu = 0.55$, $\nu = 10^{-4}$, $L = 1$, $U = 1$ and $\text{Re} = \frac{UL}{\nu}$. We carry out the simulation of NSE before turning on the $1/2$-equation model at $t^* = 1$. 

\textbf{Initial and boundary conditions.} The following initialization strategy from \cite{layton2018urans} is used:
\begin{equation*}
	k\left(x,1\right) = \frac{1}{2 \tau ^2}l^2 \left( x \right) \quad \text{where} \quad l \left( x \right) = \min \left\{0.41y, 0.082 Re ^{-1/2}\right\}
\end{equation*}
where $y$ is the wall-normal distance. According to the derivation of the $1/2$-equation model, we choose:
\begin{equation*}
	k(1) = \frac{1}{|\Omega|}\frac{1}{2\tau ^2} \int_{\Omega} l\left(x\right)^2 \, dx.\\
\end{equation*}
The turbulent viscosity $\nu_T$ is zero for $t <1$ and $t \geq 1$ is: 
\begin{equation*}
	\nu_T = \sqrt{2} \mu k(t) (\kappa y/L)^2 \tau , \quad \kappa = 0.41.
\end{equation*}
We use the BE time discretization for the momentum equation. We use the Taylor-Hood $\left( P2- P1\right)$ finite element pair for approximating the velocity field and pressure. The unstructured meshes are generated with GMSH \cite{geuzaine2009gmsh}, with a target mesh size parameter $lc$.

{\bf Order of accuracy in time.} We set target mesh size $lc = 1/36 $. Choose a very small $dt = 0.001$ to provide an approximation taken to be the true solution. The successive time steps are $dt = 0.002, 0.004, 0.006$, and $0.008$. We calculate the rate with the data from $t=1$ to $t=1.3$.
\begin{table}[!htbp]
	\caption{Errors and convergence rates in time.}
	\label{tab:1}
	\centering
	\begin{tabular}{c | c | c | c | c }
		dt & $\max_{t_n}\|u-u_{h}\|_{2,0}$ & rate & $\int_{0}^T\|\nabla u-\nabla u_{h}\|_{2,0}^2$ & rate \\
		\hline
		$8e-3$ & $0.011871$ & --& 23.53 & --\\
		$6e-3$ & $0.00897$ & 0.97 &18.50& 0.92\\
		$4e-3$  & $0.00578$ &1.08& 11.95&1.03\\
		$2e-3$ & $0.00213$ & 1.43 & 3.43 & 1.40 \\ 
	\end{tabular}
	
\end{table}

{\bf Order of accuracy in space.} We look at the ratios of differences between $\Tilde{u}_h$ computed for different $h$. We compare the solutions for the grid sizes, e.g. $h, \alpha h, \alpha^2 h$ gives
\begin{equation*}
	\frac{\Tilde{u}_h - \Tilde{u}_{\alpha h}}{\Tilde{u}_{\alpha h} - \Tilde{u}_{\alpha^2 h}} = \alpha^{-p} + O(h).
\end{equation*}
where $p$ is the order of the method \cite{OlodRunborg2012}. In our test, we take $\alpha = 3/4$. We set $dt = 0.005$ for all simulations, $h =1/60, 1/60 \cdot (3/4),/60 \cdot (3/4)^2, 1/60 \cdot (3/4)^3, 1/60 \cdot (3/4)^4$. We calculate the rates with the data from $t=1$ to $t=1.5$. 
\begin{table}[!htbp]
	\caption{Errors and convergence rates for velocity in space.}
	\label{tab:1}
	\centering
	\begin{tabular}{c | c | c | c | c }
		h & $\max_{t_n}\|u_{h}-u_{3/4\cdot h}\|_{2,0}$ & rate & $\int_{0}^T\|\nabla u_{h}-\nabla u_{3/4\cdot h}\|_{2,0}^2$ & rate \\
		\hline
		$1/60 $ & 0.045145 & -- & 2725.76 & -- \\ 
		$(\frac{3}{4})^1\cdot 1/60$  & 0.028904 &1.55 & 1922.15 &0.61\\
		$(\frac{3}{4})^2 \cdot 1/60$ & 0.011953 &3.07& 648.01 &1.89\\
		$(\frac{3}{4})^3 \cdot 1/60$ & 0.006583 & 2.07 & 230.90 & 1.79\\
	\end{tabular}
\end{table}
Rates are jumping since we are approximating a flow with non-smooth solutions and the longer the flow evolves the smaller $h$ and $\Delta t$ need to be to be in the asymptotic regime. From the above two tables, we observe the first order of accuracy in time and in average second order in space, which verifies our theoretical results.

\section{Conclusions}
Limited computation evidence in \cite{fang20231} indicates that volume averaged statistics predicted by 1-equation URANS models can be well approximated from the 1/2 equation model. This reduces computational costs provided the coupled system can be reliably and accurately approximated. We show herein that this is possible by giving a complete convergence analysis of a fundamental method and delineating how to treat the eddy viscosity nonlinearity in the numerical analysis.

\section*{Acknowledgments}
The author Wei-wei Han was partially supported by the Innovative Leading Talents Scholarship established by Xi'an Jiaotong University. The research of W. Layton and Rui Fang was partially supported by the NSF under grant DMS 2110379 and by the University of Pittsburgh Center for Research Computing through the resources provided on the SMP cluster.

\bibliography{mybib}

\end{sloppypar}
\end{document}